\documentclass[a4paper,11pt]{article}
\usepackage[T1]{fontenc}
\usepackage[english]{babel}
\usepackage{amsfonts}
\usepackage{amsmath}
\usepackage{mathabx}
\usepackage{amsthm}
\usepackage{amssymb}
\usepackage{graphicx}
\usepackage{fancyhdr}

\theoremstyle{plain}

\newtheorem{prop}{Proposition}[section]

\theoremstyle{definition}

\newtheorem{rmk}{Remark}[section]

\newcommand{\R}{\mathbb{R}}

\newcommand{\curl}{\mbox{curl}\;}
\newcommand{\dive}{\mbox{div}\;}
\newcommand{\divh}{\mbox{div}_h\;}

\newcommand{\Rey}{{\textrm{Re} \,}}
\newcommand{\Rm}{{\textrm{Rm} \,}}
\newcommand{\Ha}{\textrm{Ha} \,}
\newcommand{\pa}{\partial}
\newcommand{\na}{\nabla}
\newcommand{\bu}{\mathbf{u}}
\newcommand{\bb}{\mathbf{b}}
\newcommand{\bx}{\mathbf{x}}
\newcommand{\be}{\mathbf{e}}
\newcommand{\T}{\mathbb{T}}
\usepackage{array}
\title{Formal Derivation and Stability Analysis of \\ Boundary Layer Models in MHD}
\author{D. G\'erard-Varet, M. Prestipino}
\begin{document}
\maketitle

\begin{abstract} 
We provide a systematic derivation of boundary layer models in magnetohydrodynamics (MHD), through an asymptotic analysis of the incompressible MHD system. We recover classical linear models, related to the famous Hartmann and Shercliff layers,  as well as nonlinear ones, that we call magnetic Prandtl models. We  perform their linear stability analysis, emphasizing the stabilizing effect of the magnetic field.   
\end{abstract}

\section{Introduction}
The dynamics of an electrically conducting liquid near a wall has been a topic of constant interest, at least  since the pioneering work of Hartmann  \cite{Hart}. It is relevant to many domains of active research, such as dynamo theory \cite{DDG} or nuclear fusion \cite{Wes}. 

An appropriate starting point to describe such dynamics is the classical incompressible MHD system. It is set in an open subset $\Omega$ of $\R^3$, modeling the fluid domain.  It reads in dimensionless form   \cite{Davidson,GeLeLe}: 
\begin{equation} \label{MHD1} 
\left\{
\begin{aligned}
& \pa_t \bu  + \bu \cdot \na \bu + \na p - \frac{1}{\Rey} \Delta \bu  = S \bb \cdot \na \bb,  \\ 
& \pa_t \bb  - \curl ( \bu \times \bb) + \frac{1}{\Rm} \curl \curl  \bb  = 0,  \\
& \dive  \bu = 0, \dive \bb = 0, \quad t > 0, \quad \bx \in \Omega. 
\end{aligned}
\right.
\end{equation}
The parameters $\Rey$ and $\Rm$ are the hydrodynamic and magnetic Reynolds numbers respectively. The parameter $S$ is the so-called coupling parameter. It is given by 
$$ S = \frac{B_0^2}{\mu \rho U^2} = \frac{\Ha^2}{\Rey \Rm}, \quad \text{where } \Ha =  B_0 L \left(\frac{\sigma}{\eta}\right)^{1/2} $$
is the Hartmann number. Here, $B_0$ and $U$ are typical amplitudes for the magnetic and velocity  fields, $L$ is a typical length scale of the flow, $\rho$ is the density of the fluid,  $\mu$ is its magnetic permeability and $\eta$ is the viscosity coefficient. 

Equations in $\Omega^c$ and boundary conditions at the interface $\partial \Omega$ depend on the electrical properties of the surrounding medium $\Omega^c$. We focus here on the case of an insulator, so that
\begin{equation} \label{MHD2} 
 \curl \bb  = 0, \quad \dive \bb = 0 \quad \text{in  } \Omega^c. 
 \end{equation}
The boundary conditions at $\pa \Omega$ are  
\begin{equation} \label{BC} 
 \bu = 0, \quad [\bb]= 0 \quad \text{at } \pa \Omega.
\end{equation}
where the bracket refers to the jump of $\bb$ across the boundary $\pa \Omega$ (see \cite{Gilbert} for more).  

\medskip
For simplicity, we assume a uniform background magnetic field, meaning that $\bb = \be$ in $\Omega^c$ for some constant vector $\be$. This relation is satisfied for all times if it is satisfied initially. Under this assumption, the MHD system can be recast in $\Omega$ only: 
\begin{equation} \label{MHD}
\left\{
\begin{aligned}
& \pa_t  \bu  +  \bu \cdot \na  \bu + \na p - \frac{1}{\Rey} \Delta  \bu  = S \bb \cdot \na \bb, \\ 
& \pa_t \bb  - \curl ( \bu \times b)  + \frac{1}{\Rm} \curl \curl  \bb  = 0, \\
& \dive  \bu = 0, \dive \bb = 0, \\
&  \bu\vert_{\pa \Omega} = 0, \quad \bb\vert_{\pa \Omega} = \be.   
\end{aligned}
\right.
\end{equation}

Many MHD flows are characterized by a large hydrodynamic Reynolds number, $\Rey \gg 1$. It generates a {\em boundary layer} near $\pa \Omega$, that is a thin zone of high velocity gradients. The understanding of the boundary layer is a major problem in hydrodynamics, notably in relation to drag computation, or vortex generation. For purely hydrodynamic flows ($S = 0$ in \eqref{MHD}), a classical model for the boundary layer is the celebrated Prandtl system \cite{Prandtl}. However, this model is known to be highly unstable. It is especially true in the presence of an adverse pressure gradient, where reverse flow and boundary  layer separation can occur. 

It is then very natural to investigate the effect of a magnetic field on such instabilities. The existing results on this issue go both ways:   
\begin{itemize}
\item On one hand, stabilizing effects were stressed out. For instance, in the context of ideal MHD and plane  parallel flows, the action of a parallel magnetic field tightens the region of possible unstable wave speeds \cite{HuTo}. Another more mathematical  example is the well-posedness of inviscid hydrostatic equations between two planes, that is restored under the action of a parallel magnetic field \cite{Ren}. As regards dissipative MHD, similar stability results are known. For instance, in the regime $\Ha \gg 1$, transverse magnetic fields generate boundary layers of Hartmann type, which behave much better than the Prandtl ones \cite{Hart, Albou, Rou}. 
\item On the other hand, it was shown that magnetic fields can favour the appearance of inflexion points in the velocity profile \cite{Hunt, Nunez}. By this loss of concavity, they may generate instabilities, and one could expect earlier separation in the boundary layers. 
\end{itemize}

The purpose of this note is to gain some insight into the analysis of MHD boundary layer models. It is primarily intended to mathematicians, either applied or interested in the theory of fluid PDE's. The goal is twofold. First, we wish to  provide a clear picture of the various models available, depending on the asymptotics under consideration, and the orientation of the background field with respect to the wall.  Then, we wish to emphasize the stabilizing effect of the magnetic field, through partial linear stability analysis. We hope that this work wil serve as a starting point for more complete mathematical and numerical analysis. 

The outline of the paper is as follows. We consider the case of a half-space $\Omega = \R^3_+$, and consider both the case of a transverse and tangent background magnetic fields: $\be= \be_z$ and $\be = \be_x$, with $\bx = (x,y,z)$. The first part of the paper is a systematic derivation of MHD boundary layer models, depending on the relative scalings  of $\Rey$, $\Rm$ and $S$. We obtain in this way different sets of equations. They include linear systems, related to the classical Hartmann and Shercliff layers, but also nonlinear ones, that we call {\em magnetic Prandtl models}. 

Such magnetic Prandtl models marry features of the Prandtl equations and the Hartmann/Shercliff ones. They are interesting mathematically, because their well-posedness is unclear. Indeed, contrary to Navier-Stokes, such asymptotic models do not retain tangential diffusion. Therefore, the control of high tangential frequencies is an issue. Note that this difficulty already occurs in the classical Prandtl equation, whose well-posedness properties have been satisfactorily understood only recently \cite{GeDo, GeNg,MaWo,AlWaXuTa,GeMa,Tong1,Chen}.
In particular, for general smooth initial data, without monotonicity assumption, local well-posedness fails: it only holds under Gevrey regularity in $x$ of the data, that is under strong localization in frequency. 

In light of these results, we discuss in the second part of the paper the well-posedness of the magnetic Prandtl models. Namely, we study linearizations around shear flows, and their stability with respect to high frequencies. {\em We notably show that for tangential magnetic fields, linearizations around non-monotonic shear flows are well-posed in Sobolev spaces}. This is in sharp contrast with the Prandtl equation, which is known to be ill-posed in Sobolev spaces. Hence, while tangential magnetic fields create inflexion points in the  velocity profiles, as advocated in \cite{Hunt, Nunez}, they may at the same time suppress hydrodynamic instabilities.



 \section{Derivation of MHD layers}
 We wish to study  solutions of \eqref{MHD} that are of boundary layer type, and to find  which reduced models  they satisfy, depending on the relative values  of parameters $\Rey$, $\Rm$ and $S$. Obviously, we always assume that $\Rey \gg 1$, which is  necessary for the generation of a boundary layer. In the case $S=0$, that is in the purely hydrodynamic regime, it is well-known that a formal asymptotics leads to the so-called Prandtl equation. But of course, our goal here is rather to  emphasize the role of magnetic effects in the boundary layer:  we are interested in  models that couple equations on $\bu$ and $\bb$. Let us also stress that in most applications, the magnetic Reynolds number is usually smaller than the hydrodynamic one, so that we always assume: 
\begin{equation} \label{Reynolds_constraints}
\Rey \gg 1, \quad \Rm \lesssim \Rey. 
\end{equation} 
For simplicity, we further restrict to a simple geometry, namely the half-space $\Omega = \{ z > 0\}$. Nevertheless, we believe that our analysis could extend to curved bondaries  (through the introduction of curvilinear and transverse coordinates near the boundary).  We distinguish between the case of a transverse  
background magnetic field $\be = \be_z$ and a tangent background magnetic field, say $\be = \be_x$. 

\subsection{Layers under a transverse magnetic field} 
We consider here solutions of \eqref{MHD} behaving like: 
 \begin{equation} \label{BL}
 \begin{aligned}
 \bu & \: \approx \: \Bigl(u'_x\left(t,x,y,\lambda^{-1} z\right), \, u'_y\left(t,x,y,\lambda^{-1} z\right),  \, \lambda \, u'_z\left(t,x,y,\lambda^{-1} z\right)\Bigr), \\
 \bb & \: \approx \: \be  + \delta \Bigl(b'_x\left(t,x,y,\lambda^{-1} z\right),  \, b'_y\left(t,x,y,\lambda^{-1} z\right), \, \lambda  \, b'_z\left(t,x,y,\lambda^{-1} z\right)\Bigr)
\end{aligned}
 \end{equation}
and similarly for the pressure. The parameter $\lambda \ll 1$ denotes the size of the boundary layer: the profiles $\bu' =  (\bu'_h, u'_z) =  (u'_x, u'_y,u'_z)$, $p'$  and $\bb' = (\bb'_h, b'_z)  = (b'_x,b'_y,b'_z)$ depend on a rescaled variable 
$z' = \lambda^{-1} z$. The parameter $\delta = O(1)$ denotes the typical norm of the magnetic perturbation. Note the rescaling of the vertical components by a factor $\lambda$: it is consistent with the divergence-free conditions on $u$ and $b$.  
 
\medskip
We insert the expressions \eqref{BL} into \eqref{MHD}.  After dropping  the primes, we get 
\begin{equation} \label{MHD_rescaled} 
\left\{
\begin{aligned}
& \pa_t \bu_h  + \bu \cdot \na \bu_h + \na_h  p   - \frac{1}{\Rey} \left(\Delta_h  + \lambda^{-2} \pa^2_z\right) \bu_h  = \frac{S \delta}{\lambda} \pa_z  \bb_h  +  S \delta^2 \bb \cdot \na \bb_h, \\
& \pa_t u_z  + \bu \cdot \na u_z + \lambda^{-2}  \pa_z p   - \frac{1}{\Rey} \left(\Delta_h  + \lambda^{-2} \pa^2_z\right) u_z  = \frac{S \delta}{\lambda} \pa_z  b_z  +  S \delta^2 \bb \cdot \na b_z, \\
&  \pa_t \bb_h   - (\delta \lambda)^{-1} \pa_z \bu_h - \left( \curl (\bu \times \bb) \right)_h + \frac{1}{\Rm} \na_h \dive \bb - \frac{1}{\Rm}\left(\Delta_h + \lambda^{-2} \pa^2_z\right)  \bb_h  = 0, \\
& \pa_t b_z  - (\delta \lambda)^{-1} \pa_z u_z - \left( \curl (\bu \times \bb) \right)_z + \frac{1}{\Rm\lambda^2} \pa_z  \dive \bb  - \frac{1}{\Rm}\left(\Delta_h + \lambda^{-2} \pa^2_z\right)  b_z  = 0, \\
& \dive \bu = \dive \bb = 0,
\end{aligned}
\right.
\end{equation}
where the substrict $h$ above refers to horizontal components or variables: 
$$ \mathbf{f}_h = (f_x, f_y),  \quad \na_h = \left( \pa_x, \pa_y \right), \quad \Delta_h = \pa^2_x + \pa^2_y. $$
 The equations are completed by the Dirichlet conditions 
\begin{equation} \label{Di1}
 \bu = \bb = 0 \text{ at } z=0. 
 \end{equation}
 Moreover, we expect  vertical variations of the boundary layer solutions to be localized near $z=0$. Therefore, we impose that $\bu_h$ and  $\bb_h$ have a limit as $z \rightarrow +\infty$. 
We  denote by $\bu^\infty_h$ and $\bb^\infty_h$ such limits. We also impose that the $z$ derivatives of $\bu_h$  and $\bb_h$  decay to zero: 
\begin{equation} \label{Di2}
(\bu_h, \bb_h) \rightarrow (\bu^\infty_h, \bb^\infty_h),  \quad \pa^k_z (\bu_h, \bb_h) \rightarrow  (0,0), \quad \forall k \ge 1, \quad \text{as } z \rightarrow +\infty.
\end{equation} 
 Note that once $\bu_h$ and $\bb_h$ are determined, the divergence-free conditions and Dirichlet conditions \eqref{Di1} fully determine $u_z$ and $b_z$. From the condition \eqref{Di2}, they should be at most $O(z)$ at infinity. Note that equivalently, $b_z$ can be determined by equation (\ref{MHD_rescaled}d). This follows easily  from the well-known fact that the divergence-free condition is preserved by the evolution equations (\ref{MHD_rescaled}c,d). Indeed, taking  the divergence of (\ref{MHD_rescaled}c,d), we get 
$\pa_t \dive \bb  = 0$ in $\Omega$. 

\paragraph{Hartmann regime.}
The first case is when $\frac{S \delta}{\lambda} \gg 1$. Then, the term  $\frac{S \delta}{\lambda} \pa_z  \bb_h$ in (\ref{MHD_rescaled}a)  is diverging. It must be balanced by the term coming from diffusion in $z$. We must also keep {\it a priori} the horizontal pressure gradient, whose amplitude in the layer is unknown.  Retaining  these leading order terms, we get  
\begin{equation} \label{reducedfirsteq}
\na_h p  - \frac{1}{\Rey} \lambda^{-2} \pa^2_z \bu_h  =  \frac{S \delta}{\lambda} \pa_z  \bb_h,
\end{equation}
which yields in particular that 
\begin{equation} \label{scaling1}
 \frac{1}{\Rey} \lambda^{-2} \sim \frac{S \delta}{\lambda}.  
 \end{equation}
With this balance and the assumption $\Rey \gg 1$, the second equation (\ref{MHD_rescaled}b)
  yields at leading order: 
$\pa_z p = 0$. We recover the classical fact that the pressure is constant in boundary layers.  Back to \eqref{reducedfirsteq}, we can send $z$ to infinity and use \eqref{Di2} to deduce that $\na_h p = 0$ and 
\begin{equation} \label{reducedfirsteq_2}
- \frac{1}{\Rey} \lambda^{-2} \pa^2_z \bu_h  =  \frac{S \delta}{\lambda} \pa_z  \bb_h,
\end{equation}
Similarly, in (\ref{MHD_rescaled}c), the only term that can balance $\lambda^{-1} \pa_z \bu_h$ is the term coming from diffusion in $z$. Retaining these two terms we get 
\begin{equation} \label{reducedthirdeq}
  - (\delta \lambda)^{-1} \pa_z \bu_h  - \frac{1}{\Rm} \lambda^{-2} \pa^2_z \bb_h  = 0, 
\end{equation}
so that 
\begin{equation} \label{scaling2}
 (\delta \lambda)^{-1}  \sim \frac{1}{\Rm} \lambda^{-2}. 
 \end{equation}
Combining \eqref{scaling1} and \eqref{scaling2}, we get 
$$ \lambda^2 \sim \frac{1}{\Rey \, \Rm S} \sim \Ha^{-2}.  $$
Hence, {\em the typical size of the layer is $\Ha^{-1}$}. We set 
\begin{equation} 
\lambda = \Ha^{-1}, \quad \delta = \Rm \Ha^{-1}. 
\end{equation} 
The previous equations \eqref{reducedfirsteq_2}-\eqref{reducedthirdeq} on $\bu_h,\bb_h$ simplify into 
\begin{equation} \label{Hartmann}
\pa^2_z \bu_h + \pa_z \bb_h = 0,  \quad  \pa_z \bu_h  +  \pa^2_z \bb_h = 0 
\end{equation}    
which yields 
$$ - \pa^3_z \bu_h = - \pa_z \bu_h. $$   
From the boundary conditions, we deduce
\begin{equation} \label{Hartmannbis}
\bu_h \: = \:  (1 - e^{-z}) \bu^\infty_h, \quad  \bb_h = (1- e^{-z}) \bu^\infty_h 
\end{equation} 
or
$$ \bu_h \: = \:  (1 - e^{- \Ha z}) \bu^\infty_h, \quad  \bb_h = (1- e^{-\Ha z}) \bu^\infty_h.  $$
in the original $z$ variable.  These are the classical Hartmann profiles.  
\begin{rmk}
As the focus of our note is on boundary layers, we do not adress the dynamics of the limits at infinity $\bu^\infty_h(t,x,y)$ and $\bb^\infty_h(t,x,y)$. In  a full analysis of \eqref{MHD}, these limits appear as the boundary values  of velocity  and magnetic   fields  $\bu_h^{int}$ and  $\bb_h^{int}$, describing the (horizontal) dynamics away from the boundary layer. Hence, they are not arbitrary, but constrained by equations \eqref{MHD} and the solvability of the boundary layer. For instance, in the Hartmann case, we see from \eqref{Hartmannbis} that one condition is  $\bb^\infty_h = \bu^\infty_h$. 
\end{rmk}

\begin{rmk}
To be consistent, the derivation of the Hartmann boundary layer requires {\it a priori} some  assumptions on the parameters. The first requirement is of course that the size $\lambda$ of the layer be small, or equivalently $\Ha \gg 1$. Also, we assumed that  $\frac{S \delta}{\lambda} \gg 1$, that is $\frac{\Ha^2}{\Rey} \gg 1$. Eventually, the condition $\delta = O(1)$ means $\Rm \Ha^{-1} = O(1)$. Note however that this last condition on $\delta$ is not needed in the derivation of the Hartmann equations: a sufficient condition is that  $\frac{S\delta}{\lambda} \gg S \delta^2$ and $(\delta \lambda)^{-1} \gg 1$. Both conditions come down to 
$\Ha^2 \gg \Rm$,  which is automatically satisfied if $\Ha^2 \gg \Rey$ and  $\Rm  \lesssim \Rey$ (see \eqref{Reynolds_constraints}). Note also that these assumptions can  be sometimes relaxed. For instance, in the case where $\bu^\infty_h$ is constant, one can check that the Hartmann profiles \eqref{Hartmannbis} are exact solutions of the full system \eqref{MHD_rescaled} (with $u_z = b_z = 0$). 
 \end{rmk}
 
 \paragraph{Mixed Prandtl/Hartmann regime.} The second case is when $\frac{S \delta}{\lambda} \sim 1$. In this case, the convective term in the equation for $\bu_h$ can no longer be neglected. Hence, the leading order dynamics reads:
 \begin{equation} \label{Prandtl_MHD}
   \pa_t \bu_h  + \bu \cdot \na \bu_h + \na_h  p   - \frac{1}{\Rey}  \lambda^{-2} \pa^2_z \bu_h  = \frac{S \delta}{\lambda} \pa_z  \bb_h, 
   \end{equation}
Meanwhile,  the induction equation still yields  the same balance:  
$$ - (\delta \lambda)^{-1} \pa_z \bu_h  - \frac{1}{\Rm} \lambda^{-2} \pa^2_z \bb_h  = 0, $$
  or after integration in $z$: 
\begin{equation} \label{reduced_induction}
 - (\delta \lambda)^{-1}  \bu_h  - \frac{1}{\Rm} \lambda^{-2} \pa_z \bb_h  = 0. 
 \end{equation}
As before, we can  take  $ \lambda = \Ha^{-1}$. Note  that $\frac{1}{\Rey} \lambda^{-2} \sim \frac{S \delta}{\lambda}  \sim 1$,  giving the extra condition 
$$\lambda^{2} \sim \Rey^{-1},  \text{ or }   \Ha \sim \sqrt{\Rey}. $$  
 Moreover, the  equation for the vertical velocity component gives at leading order: $ \pa_z p = 0$. Eventually, substituting \eqref{reduced_induction} in \eqref{Prandtl_MHD}, we obtain the system
 \begin{equation} \label{dampedP}
 \left\{
 \begin{aligned}
& \pa_t \bu_h  + u \cdot \na \bu_h + \na_h  p   - \frac{\Ha^2}{\Rey}  \pa^2_z \bu_h  + \frac{\Ha^2}{\Rey}  \bu_h = 0, \\ 
 & \pa_z p = 0, \\ 
 & \divh \bu_h  + \pa_z u_z = 0. 
 \end{aligned}
 \right. 
 \end{equation} 
 We recognize a nonlinear Prandtl type equation, {\em with an extra magnetic damping term}. This model belongs to what we called in the introduction {\em magnetic Prandtl models}, mixing features of Prandtl and Hartmann dynamics.

 \subsection{Layers in a tangent magnetic field.}  
 In this section, we consider the case of a tangent background magnetic field $\bb = \be_x$. As the MHD system is invariant through horizontal rotation, the choice of $\be_x$ is no loss of generality. Proceeding as before, we look for approximate solutions of the type
  \begin{equation} \label{BL2}
 \begin{aligned}
 \bu & \: \approx \: \Bigl(u'_x\left(t,x,y,\lambda^{-1} z\right), \, u'_y\left(t,x,y,\lambda^{-1} z\right),  \, \lambda \, u'_z\left(t,x,y,\lambda^{-1} z\right)\Bigr), \\
 \bb & \: \approx \: \be_x + \delta \Bigl(b'_x\left(t,x,y,\lambda^{-1} z\right),  \, b'_y\left(t,x,y,\lambda^{-1} z\right), \, \lambda  \, b'_z\left(t,x,y,\lambda^{-1} z\right)\Bigr).  
\end{aligned}
 \end{equation}
By plugging these approximations in the MHD equations, we have this time: 
\begin{equation} \label{MHD_rescaled2} 
\left\{
\begin{aligned}
& \pa_t \bu_h  + \bu \cdot \na \bu_h + \na_h  p   - \frac{1}{\Rey} \left(\Delta_h  + \lambda^{-2} \pa^2_z\right) \bu_h  = S \delta \pa_x  \bb_h  +  S \delta^2 \bb \cdot \na \bb_h, \\
& \pa_t u_z  + \bu \cdot \na u_z + \lambda^{-2}  \pa_z p   - \frac{1}{\Rey} \left(\Delta_h  + \lambda^{-2} \pa^2_z\right) u_z  = S \delta \pa_x  b_z  +  S \delta^2 \bb \cdot \na b_z, \\
&  \pa_t \bb_h   - \delta^{-1} \pa_x \bu_h - \left( \curl (\bu \times \bb) \right)_h + \frac{1}{\Rm} \na_h \dive \bb - \frac{1}{\Rm}\left(\Delta_h + \lambda^{-2} \pa^2_z\right)  \bb_h  = 0, \\
& \pa_t b_z  - \delta^{-1} \pa_x u_z - \left( \curl (\bu \times \bb) \right)_z + \frac{1}{\Rm \lambda^2} \pa_z  \dive \bb  - \frac{1}{\Rm}\left(\Delta_h + \lambda^{-2} \pa^2_z\right)  b_z  = 0, \\
& \dive \bu = \dive \bb = 0.
\end{aligned}
\right.
\end{equation}
This system is still completed by \eqref{Di1}-\eqref{Di2}. Note that when $\delta \sim 1$, the last two terms at the right-hand side of (\ref{MHD_rescaled2}a,b) have the same amplitude. The same remark applies to the terms 
$\delta^{-1} \pa_x \bu$  and $\curl (\bu \times \bb)$, see the third and fourth equations. In other words, when $\delta \sim 1$, the perturbative writing (\ref{BL2}b) is somehow artificial, and should be replaced by 
$$  \bb  \: \approx \:  \Bigl(b'_x\left(t,x,y,\lambda^{-1} z\right),  \, b'_y\left(t,x,y,\lambda^{-1} z\right), \, \lambda  \, b'_z\left(t,x,y,\lambda^{-1} z\right)\Bigr).  $$
 We shall consider this non perturbative regime at the end of the section. 
 
 \paragraph{Shercliff regime.}
 We consider here that 
 $$ \delta \ll 1, \quad S \delta  \gg 1. $$
 In the equation for $\bu_h$, the diffusion in $z$ and the horizontal pressure gradient can balance the linearized Lorentz force $S \delta \pa_x  \bb_h$. The reduced dynamics reads
\begin{equation} \label{reducedfirsteqbis}
 \na_h p - \frac{1}{\Rey} \lambda^{-2} \pa^2_z \bu_h  = S \delta \pa_x  \bb_h  
\end{equation}
 and in particular
 \begin{equation} \label{scaling3}
  \frac{1}{\Rey} \lambda^{-2} \sim S \delta. 
  \end{equation}
Like in the Hartmann regime, the second equation yields at leading order $\pa_z p = 0$. Taking into account \eqref{Di2}, we  then rewrite equation \eqref{reducedfirsteqbis} as 
\begin{equation} 
S \delta \pa_x \bb_h^\infty  - \frac{1}{\Rey} \lambda^{-2} \pa^2_z \bu_h  = S \delta \pa_x  \bb_h  
\end{equation}
Similarly,  in the equation for $\bb_h$, only the magnetic diffusion in $z$ can balance $- \delta^{-1} \pa_x \bu_h$. We find 
$$ \delta^{-1} \pa_x \bu_h  + \frac{1}{\Rm} \lambda^{-2} \pa^2_z  \bb_h = 0$$
and in particular
 \begin{equation} \label{scaling4}
\delta^{-1} \sim  \frac{1}{\Rm} \lambda^{-2}  
\end{equation} 
Combining \eqref{scaling3} and \eqref{scaling4} yields $\lambda^4 \sim \Ha^{-2}$. We set 
$$ \lambda = \Ha^{-1/2}. $$
The previous equations resume to 
\begin{equation} \label{shercliff}
  \pa_x (\bb_h - \bb_h^\infty) + \pa^2_z \bu_h  = 0, \quad  \pa_x \bu_h  + \pa^2_z  \bb_h = 0. 
 \end{equation} 
 These equations describe  the so-called Shercliff layer, of typical size $\Ha^{-1/2}$ \cite{Shercliff}. In the half-space case, they can be solved by taking the Fourier transform in $x$. Accounting for \eqref{Di1}-\eqref{Di2}, we find
\begin{align*}
 \widehat{\bu_h}(\xi,z) & = - i \frac{\xi}{|\xi|}  \widehat{\bb_h^\infty}(\xi) e^{-\sqrt{\frac{|\xi|}{2}} z} \sin\left(\sqrt{\frac{|\xi|}{2}} z\right)  \\
 \widehat{\bb_h}(\xi,z) & =  \widehat{\bb_h^\infty}(\xi)  \left( 1 -  e^{-\sqrt{\frac{|\xi|}{2}} z} \cos\left(\sqrt{\frac{|\xi|}{2}} z\right) \right). 
 \end{align*}

 \begin{rmk}
 In this derivation, we assumed implicitly that $\lambda \ll 1$, that is  $\Ha \gg 1$. Also, we assumed that $\delta \ll 1$, which amounts to $\Rm \Ha^{-1} \ll 1$, as well as $S \delta \gg 1$, which amounts to $\Ha \gg \Rey$. Taking \eqref{Reynolds_constraints} into account, the  constraint $\Ha \gg \Rey$ is the more stringent. 
 \end{rmk}
 \paragraph{Mixed Prandtl/Shercliff regime.} 
We still assume here that $\delta \ll 1$, but  $S \delta  \sim 1$. One must then retain all terms of order one  in the equation for $\bu_h$, namely 
$$   \pa_t \bu_h  + \bu \cdot \na \bu_h + \na_h  p   - \frac{1}{\Rey} \lambda^{-2} \pa^2_z  \bu_h  = S \delta \pa_x  \bb_h. $$
The leading order terms in the equation for $\bb_h$ remain the same: 
$$ \delta^{-1} \pa_x \bu_h  + \frac{1}{\Rm} \lambda^{-2} \pa^2_z  \bb_h = 0.$$
It is therefore legitimate to maintain the same definition for the boundary layer size, that is 
$\lambda = \Ha^{-1/2}$. As $ \frac{1}{\Rey} \lambda^{-2}  \sim   S \delta  \sim 1$, the regime that we investigate here corresponds to 
$$\Rey \sim \Ha. $$ 
 We finally obtain the following boundary layer system: 
 \begin{equation}  \label{mixedSP}
 \left\{
 \begin{aligned}
 & \pa_t \bu_h  + \bu \cdot \na \bu_h + \na_h  p   - \frac{\Ha}{\Rey}  \pa^2_z  \bu_h  = \frac{\Ha}{\Rey}  \pa_x  \bb_h, \\
 & \pa_z  p = 0, \\
 & \pa_x \bu_h  +\pa^2_z  \bb_h = 0, \\
 & \dive \bu = 0. 
\end{aligned}
 \right. 
\end{equation}
This is a mixed Prandtl/Shercliff system. 

\paragraph{Fully nonlinear MHD layer.}
We eventually consider the case where the perturbation to the constant magnetic field $\be_x$ is of size one. In such setting, distinguishing between $\be_x$ and its perturbation is artificial. One rather looks directly for 
$$  \bb  \: \approx \:   \Bigl(b'_x\left(t,x,y,\lambda^{-1} z\right),  \, b'_y\left(t,x,y,\lambda^{-1} z\right), \, \lambda  \, b'_z\left(t,x,y,\lambda^{-1} z\right)\Bigr).  $$
We plug this new expansion into \eqref{MHD}, to obtain
\begin{equation} \label{MHD_rescaled3} 
\left\{
\begin{aligned}
& \pa_t \bu_h  + \bu \cdot \na \bu_h + \na_h  p   - \frac{1}{\Rey} \left(\Delta_h  + \lambda^{-2} \pa^2_z\right) \bu_h  =   S \bb \cdot \na \bb_h, \\
& \pa_t u_z  + \bu \cdot \na u_z + \lambda^{-2}  \pa_z p   - \frac{1}{\Rey} \left(\Delta_h  + \lambda^{-2} \pa^2_z\right) u_z  = S  \bb \cdot \na b_z, \\
&  \pa_t \bb_h   - \left( \curl (\bu \times \bb) \right)_h + \frac{1}{\Rm} \na_h \dive \bb - \frac{1}{\Rm}\left(\Delta_h + \lambda^{-2} \pa^2_z\right)  \bb_h  = 0, \\
& \pa_t b_z  - \left( \curl (\bu \times \bb) \right)_z + \frac{1}{\Rm \lambda^2} \pa_z  \dive \bb   - \frac{1}{\Rm}\left(\Delta_h + \lambda^{-2} \pa^2_z\right)  b_z  = 0, \\
& \dive \bu = \dive\bb = 0.
\end{aligned}
\right.
\end{equation}
We stress that the Dirichlet conditions are now 
\begin{equation} \label{Dir_nonlinear}
 \bu = 0, \quad \bb = \be_x \quad \text{ at } z=0. 
 \end{equation}
Let us first consider  the case  $S \gg 1$. On one hand, the contribution of the  Lorentz force diverges in (\ref{MHD_rescaled3}a), and is expected to be balanced by  the diffusion in $z$, resulting in  
$$ \frac{1}{\Rey} \lambda^{-2} \sim S  \gg 1.$$
On the other hand, looking at  the equation (\ref{MHD_rescaled3}c), we see that  
$$ \frac{1}{\Rm} \lambda^{-2} \lesssim 1$$  
otherwise the dynamics of $\bb_h$ would be trivial.  But the constraints 
$ \frac{1}{\Rey} \lambda^{-2}  \gg 1 \quad \text{ and }  \frac{1}{\Rm} \lambda^{-2}  \lesssim 1 $
are incompatible with \eqref{Reynolds_constraints}. 

\medskip
\noindent
The only relevant case is therefore $S \sim 1$: the case $S \ll 1$, leading to the usual Prandtl equation, does not exhibit any magnetic effect. To be consistent with the Dirichlet conditions, the reduced boundary layer model should contain diffusion terms for both the velocity and the magnetic field. This is possible under the two conditions
$$ \frac{1}{\Rey} \lambda^{-2} \sim S  \sim 1, \quad \frac{1}{\Rm} \lambda^{-2} \sim 1 $$
which imply 
$$ \Rey \sim \Rm \sim \Ha, \quad \lambda \sim \frac{1}{\sqrt{\Rey}}. $$ 
We set $\lambda = \frac{1}{\sqrt{\Rey}}$. We find the MHD boundary layer system
\begin{equation*} 
\left\{
\begin{aligned}
& \pa_t \bu_h  + \bu \cdot \na \bu_h + \na_h  p   -   \pa^2_z  \bu_h  =   S \bb \cdot \na \bb_h, \\
&   \pa_z p  = 0, \\
&  \pa_t \bb_h   - \left( \curl (\bu \times \bb) \right)_h - \frac{\Rey}{\Rm} \pa^2_z  \bb_h  = 0, \\
& \pa_t b_z  - \left( \curl (\bu \times \bb) \right)_z + \frac{\Rey}{\Rm} \pa_z  \dive \bb  - \frac{\Rey}{\Rm} \pa^2_z   b_z  = 0, \\
& \dive \bu = \dive \bb = 0.
\end{aligned}
\right. 
\end{equation*}
As discussed before, the divergence-free condition on $\bb$ is preserved by the evolution equation on $(\bb_h,b_z)$, so that we can get rid of the equation $\dive \bb = 0$ in the previous system. On the contrary, if we keep this equation, we can set the term $\frac{\Rey}{\Rm} \pa_z  \dive \bb$ to zero in the equation for $b_z$, and the MHD boundary layer system then reads 
\begin{equation} \label{MHDBL} 
\left\{
\begin{aligned}
& \pa_t \bu_h  + \bu \cdot \na \bu_h + \na_h  p   -   \pa^2_z  \bu_h  =   S \bb \cdot \na \bb_h, \\
&   \pa_z p  = 0, \\
&  \pa_t \bb   - \left( \curl (\bu \times \bb) \right) - \frac{\Rey}{\Rm} \pa^2_z  \bb  = 0, \\
& \dive \bu = \dive \bb = 0.
\end{aligned}
\right. 
\end{equation}
\begin{rmk}
The derivation of  \eqref{MHDBL} as an asymptotic boundary layer model is only valid under stringent assumptions on the coupling parameter and the Reynolds numbers: 
$$ \Rey \sim \Rm \sim \Ha \gg 1.$$
Still, compared to the two models derived earlier (the Shercliff and Prandtl/Shercliff systems), it is the one that retains most terms from the original system \eqref{MHD}. The other two can be seen as degeneracies from it.  
\end{rmk}
\subsection{Summary of the formal derivation}
To gather the results of the previous paragraphs, we draw  the following table, that relates the various boundary layer models to the various asymptotic regimes and  to the orientation of the magnetic field: 

\medskip
\begin{tabular}{| c | c | c | c |}
\cline{2-4}
\multicolumn{1}{ c |}{} & Linear models &  \multicolumn{2}{c |}{Nonlinear models} \\
\hline
Transverse field   & $\Ha^2 \gg \Rey$  & \multicolumn{2}{c |}{$\Ha^2 \sim \Rey$} \\
(Layer size $\Ha^{-1}$) & Hartmann, {\it cf} \eqref{Hartmann}. &  \multicolumn{2}{c |}{damped Prandtl, {\it cf}  \eqref{dampedP}.} \\
\hline
 Tangent  field & $\Ha \gg \Rey$ & $\Ha \sim \Rey \gg \Rm$&  $\Ha \sim \Rey \sim \Rm$ \\
 (Layer size $\Ha^{-1/2}$)& Shercliff , {\it cf}  \eqref{shercliff}. & mixed Prandtl/Shercliff, {\it cf}  \eqref{mixedSP}. & fully nonlinear, {\it cf}  \eqref{MHDBL}. \\
 \hline
 \end{tabular}
\section{Linear Stability}
The previous derivation is of course formal. It assumes the existence of solutions of  \eqref{MHD} that take the approximate form \eqref{BL} and \eqref{BL2}.  To ground this idea on rigorous arguments, two further steps are needed:
\begin{itemize}
\item To show that the reduced boundary layer models are well-posed, at least locally in time, so that boundary layer expansions can be built.
\item To show that once they are built, these expansions are good approximations of exact MHD solutions, over some reasonable time. This is a stability issue within the MHD system \eqref{MHD}. 
\end{itemize}
We shall provide here elements for the first step only. For simplicity, we will assume invariance with respect to $y$, and restrict in this way to two-dimensional boundary layer models: $x \in \T $, $z > 0$.   Let us note that for the  classical 2D Prandtl system, with velocity field $\bu = (u,v)$, 
\begin{equation} \label{P}
\begin{aligned}
 \pa_t u + u \pa_x u + v \pa_z u  - \pa^2_z u +\pa_x p & = 0, \\   
\pa_z p & = 0, \\
\pa_x u + \pa_z v & = 0, \\ 
u\vert_{z=0} = v\vert_{z=0} & = 0, \\
 u \rightarrow u^\infty, \quad p \rightarrow p^\infty & \quad \text{ as } z \rightarrow +\infty, 
\end{aligned}
\end{equation}
the well-posedness theory is already difficult, and was only recently well-understood. To explain the underlying difficulties, it is worth considering simple linearizations, say around shear flows:  $u =  U(z), v = 0$.  Linearized Prandtl  then reads   
\begin{equation} \label{LP}
\begin{aligned} 
 \pa_t u + U \pa_x u +  v U'  -  \pa^2_z u & = 0, \\
\pa_x u + \pa_z v & = 0, \\  
u\vert_{z=0} = v\vert_{z=0}  & = 0, \\
u \rightarrow 0 &  \quad \text{ as } z \rightarrow +\infty, 
\end{aligned}
\end{equation}
where $(u,v)$  now  refers the perturbation. The main  problem comes from the term  $v U'$: indeed, in the Prandtl model, $v$ is recovered from $u$ through the divergence-free condition:  $v = - \int_0^z \pa_x u$. This is a  first order term in $u$ (with respect to variable $x$), and contrary to the transport term  $U \pa_x u$ it has no hyperbolic structure. Hence, no basic energy  estimate can be achieved. Indeed, it turns out that the  $L^2$ type well-posedness of \eqref{LP} requires a monotonicity assumption on the velocity profile $U$. Let us stress that a similar monotonicity assumption is needed on the initial data for the nonlinear system \eqref{P} to be well-posed in Sobolev spaces, see for instance \cite{MaWo}. On the contrary, when $U$ has a  non-degenerate critical point $a$, system \eqref{LP} is ill-posed in  $L^2$ or Sobolev regularity:  it has solutions that behave like 
$$ u \approx e^{i k x} e^{i \omega(k) t} U_k(z), \quad \text{ with } \omega(k) = - k U(a) + \sqrt{|k|} \tau, \quad \Im \tau < 0, \quad  |k| \gg 1, $$  
see \cite{Cow,GeDo}. Hence, it admits unstable modes whose growth rate is proportional to the square root of the wave number $k$. As a consequence, the only functional settings that can be preserved by the Prandtl evolution in small time are made of functions highly localized in frequency: their Fourier mode $k$  in $x$ should decay at least like $e^{-\delta \sqrt{|k|}}$ for some $\delta > 0$. This corresponds to Gevrey 2 regularity in $x$. Accordingly, local well-posedness  results  in such Gevrey classes were obtained recently for the full Prandtl system: see \cite{GeMa,Tong1,Chen}. 

On the basis of these results in the hydrodynamic case, it is very interesting to investigate the effect of the magnetic field on boundary layer stability, and notably the well-posedness of MHD boundary layer models. Following the previous sections, we can distinguish between linear and nonlinear models. The  two linear models that we have derived are the Hartmann system \eqref{Hartmann} and the Shercliff system \eqref{shercliff}. They do not raise any mathematical difficulty. System \eqref{Hartmann} is made of ODEs in variable $z$, and can be solved explicitly. The same is true for \eqref{shercliff} after Fourier transform in variable $x$. The variable $t$ is only a parameter and appears through the functions $\bu_h$ and  $\bb_h$, that is through the dynamics outside the boundary layer.  

From the point of view of well-posedness, the interesting systems are the nonlinear ones, that mix Prandtl and magnetic features. We call them magnetic Prandtl models. They correspond to equations \eqref{dampedP} (with background transverse magnetic field $\be = \be_z$), \eqref{mixedSP} and \eqref{MHDBL} (with background tangential magnetic field $\be = \be_x$). We shall discuss their well-posedness properties in the next section. As explained above, {\em we shall restrict to the 2D case in variables $(x,z)$, with $\bu = (u,v)$, $\bb = (b,c)$}.  The 3D case could carry additional difficulties, see \cite{Tong3} in the classical Prandtl case. 

\subsection{Mixed Prandtl/Hartmann regime}
The 2D version of \eqref{dampedP} reads 
\begin{equation} \label{2DdampedP}
 \left\{
 \begin{aligned}
 \pa_t u + u \pa_x u + v \pa_z u  - \frac{\Ha^2}{\Rey} \pa^2_z u   + \frac{\Ha^2}{\Rey}  u &  = -\pa_x p^\infty , \\
  \pa_x u  + \pa_z v & = 0, \\ 
 u\vert_{z=0} = v\vert_{z=0} & = 0, \\
 u \rightarrow u^\infty & \quad \text{ as } z \rightarrow +\infty.
 \end{aligned}
 \right. 
 \end{equation} 
We recall that $u^\infty, p^\infty$ are known functions of $t$ and $x$, which are the trace of an Euler flow: they satisfy  
$$ \pa_t u^\infty + u^\infty \pa_x u^\infty = - \pa_x p^\infty. $$
The only difference with the usual Prandtl system is the additional damping  $\frac{\Ha^2}{\Rey}  u$. 

This damping does not affect the usual well-posedness theory (or  in other words the stability properties of high frequencies). A close look at papers \cite{GeMa,Tong1,GeDo} shows  that both the Gevrey well-posedness results and the Sobolev ill-posedness results apply to \eqref{2DdampedP}. 

\subsection{Mixed Prandtl/Shercliff regime}
The 2D version of \eqref{mixedSP} reads 
\begin{equation}  \label{2DmixedSP}
 \left\{
 \begin{aligned}
 & \pa_t u + u \pa_x u + v \pa_z u   - \frac{\Ha}{\Rey}  \pa^2_z u  = \frac{\Ha}{\Rey}  \pa_x b - \pa_x p^\infty, \\
 &  \pa_x u  +\pa^2_z b   = 0, \\
 & \pa_x u + \pa_z v = 0, \\
 &   u\vert_{z=0} = v\vert_{z=0}  = b\vert_{z=0} = 0, \\
 & u \rightarrow u^\infty, \quad b \rightarrow b^\infty, \quad \text{ as } z \rightarrow +\infty. 
\end{aligned}
 \right. 
\end{equation}

Contrary to the simple damping term due to a transverse magnetic field, the effect created by a tangential magnetic field is more subtle. Strikingly, in the context of \eqref{2DmixedSP}, it is stabilizing.  To provide a clear illustration of this fact, we restrict ourselves to a simple linearization, namely around 
$$ u = U(z), \quad v = 0,  \quad b = b^\infty \text{ constant}. $$
We assume that $U$ connects $0$ at $z =0$ to some constant $u^\infty$ at infinity. The linearized system reads 
\begin{equation} \label{LSP} 
\left\{
 \begin{aligned}
 \pa_t u + U \pa_x u +  v U'  -  \frac{\Ha}{\Rey} \pa^2_z u & = \frac{\Ha}{\Rey}  \pa_x b, \\
 \pa_x u  +\pa^2_z b  & = 0, \\
\pa_x u + \pa_z v & = 0, \\  
u\vert_{z=0} = v\vert_{z=0} =  b\vert_{z=0} & = 0, \: (u,b) \rightarrow  0 &   \quad \text{ as } z \rightarrow +\infty.
\end{aligned}
\right.
\end{equation}
Our aim is to  prove good {\it a priori} estimates for this linear system, in the Sobolev framework. Therefore, we introduce the analogue of vorticity, which in the boundary layer  context is simply $\omega = \pa_z u$. Differentiating the first equation with respect to $z$, we find 
$$ \pa_t \omega + U \pa_x \omega +  v U''  -  \frac{\Ha}{\Rey} \pa^2_z \omega  = \frac{\Ha}{\Rey}  \pa_x \pa_z b. $$
We remark that $\pa_z \omega\vert_{z=0} = \pa^2_z u\vert_{z=0} = 0$, as can be seen from evaluating (\ref{LSP}a) at $z=0$. Multiplication by $\omega$ and integration over $\Omega = \T \times \R_+$ give 
$$ \frac{1}{2} \frac{d}{dt} \| \omega \|_{L^2}^2 +  \frac{\Ha}{\Rey} \| \pa_z \omega \|_{L^2}^2 = - \int_{\Omega} U'' v \omega + \frac{\Ha}{\Rey} \int_{\Omega} \pa_x \pa_z b \, \omega. $$
The first term at the r.h.s. is bounded by 
$$   \left| \int_{\Omega} U'' v \omega \right| \le \| U'' \int_0^z \pa_x u  \|_{L^2} \, \|\omega \|_{L^2} \: \le \:  2 \|z U'' \|_{L^\infty}  \|\pa_x u \|_{L^2}  \| \omega \|_{L^2}, $$
where we assumed implicitly that $z \rightarrow z U''$ is bounded and  applied the Hardy inequality to the first factor. As regards the additional term, we use the second equation to get 
$$ \int_{\Omega} \pa_x \pa_z b \, \omega =  - \int_{\Omega} \pa_x \pa^2_z b  \,  u =  \int_{\Omega} \pa^2_x u \, u  =  - \int_{\Omega} | \pa_x u|^2. $$
Hence, we get 
$$  \frac{1}{2} \frac{d}{dt} \| \omega \|_{L^2}^2 +  \frac{\Ha}{\Rey} \left( \| \pa_z \omega \|^2_{L^2} + \| \pa_x u \|_{L^2}^2 \right)  \le 2 \|z U'' \|_{L^\infty} \|\pa_x u \|_{L^2}  \| \omega \|_{L^2} $$
 which implies
\begin{equation} \label{ineq_omega}
\frac{1}{2}\frac{d}{dt} \| \omega \|_{L^2}^2 +   \frac{\Ha}{2 \Rey} \left( \| \pa_z \omega \|_{L^2}^2  +   \|\pa_x u \|_{L^2}^2 \right)  \le  C \|\omega \|_{L^2}^2,
\end{equation}
with $C =   2 \|z U'' \|_{L^\infty}^2 \frac{\Rey}{\Ha}$.
To have some information on $u$ itself rather than $\omega$, we  perform another energy estimate directly on (\ref{LSP}a), which gives
$$  \frac{1}{2} \frac{d}{dt} \| u \|_{L^2}^2 +  \frac{\Ha}{\Rey} \| \pa_z u \|_{L^2}^2 = - \int_{\Omega} U' v \, u + \frac{\Ha}{\Rey} \int_{\Omega} \pa_x  b \, u. $$
As previously, we have 
 $$    \left| \int_{\Omega} U'  v  \, u \right|  \: \le \:  2 \|z U' \|_{L^\infty}  \|\pa_x u \|_{L^2}  \| u \|_{L^2}, \quad 
  \int_{\Omega} \pa_x  b \, u  = - \int_\Omega |\pa_z b|^2$$
 and we end up with 
\begin{equation} 
\frac{1}{2} \frac{d}{dt} \|u \|_{L^2}^2 +  \frac{\Ha}{\Rey} \left( \|\pa_z u \|_{L^2}^2  +   \|\pa_z b \|_{L^2}^2  \right) \le  \|z U' \|_{L^\infty} \left( \|\pa_x u \|_{L^2}^2 + \|u  \|_{L^2}^2 \right). 
\end{equation}
Combining with inequality \eqref{ineq_omega}, we obtain
$$ \frac{1}{2} \frac{d}{dt} \left( \|u \|_{L^2}^2  + (1 + \alpha) \| \omega \|_{L^2}^2 \right)  + \frac{\Ha}{2\Rey} \left(   \| \pa_z b \|_{L^2}^2  +   \| \pa_z \omega \|_{L^2}^2  +   \|\pa_x u \|_{L^2}^2 \right) \le C' \left( \| u \|_{L^2}^2  + \| \omega \|_{L^2}^2 \right),  $$
where $\alpha = \frac{2 \Rey}{\Ha} \| z U' \|_{L^\infty}$, and $C'  = \max (\|z U' \|_{L^\infty}, C(1+\alpha))$.   
Eventually, with Gronwall inequality: 
\begin{equation} \label{Gronwall1}
\|\omega(t) \|^2  + \| u(t) \|_{L^2}^2  + \int_0^t   (\| \pa_z b \|_{L^2}^2  +  \| \pa_z \omega \|_{L^2}^2  +   
 \| \pa_x u \|_{L^2}^2) \le M (\|\omega_0 \|_{L^2}^2 +  \| u_0 \|_{L^2}^2)   e^{M t}, \quad \forall t \ge 0,  
\end{equation}
where $M > 0$ is large enough.  Eventually, to have some more information on $b$, one can multiply (\ref{LSP}a) by $\pa_t u$. Integrating over $\Omega$ and over $[0,t]$, we get after straightforward manipulations: 
$$  \int_0^t  \| \pa_t u \|_{L^2}^2  - \frac{\Ha}{\Rey} \int_0^t  \int_\Omega \pa_x b \, \pa_t u  \:  \le \:  C  \int_0^t  (\| \pa_x u \|_L^2 + \| \pa_z \omega \|_{L^2} ) \| \pa_t u \|_{L^2}.  $$
Using (\ref{LSP}b), we find 
$$  \int_0^t  \int_\Omega \pa_x b \, \pa_t u =  \int_0^t \int_\Omega b \, \pa_t \pa^2_z b = \frac{1}{2} \| \pa_z b(t) \|_{L^2}^2 - \frac{1}{2} \| \pa_z b_0 \|_{L^2}^2,  $$
and can conclude that 
\begin{equation} \label{Gronwall2}
 \int_0^t \| \pa_t u \|_{L^2}^2  +   \frac{\Ha}{\Rey}  \| \pa_z b(t) \|_{L^2}^2  \le   \frac{\Ha}{\Rey} \|\pa_z b_0 \|_{L^2}^2 +  C^2 M (\|\omega_0 \|_{L^2}^2 +  \| u_0 \|_{L^2}^2)  e^{M t}, \quad \forall t \ge 0.
  \end{equation}
 Let us stress that, from the bounds \eqref{Gronwall1} and \eqref{Gronwall2}, all terms at the l.h.s. of (\ref{LSP}a) belong to $L^2_{loc}(\R_+, L^2(\Omega))$, and therefore so does the r.h.s. $\pa_x b$. Moreover, $\pa_z b$ belongs to $L^\infty_{loc}(\R_+, L^2(\Omega))$, as seen from \eqref{Gronwall2}.  We recall that $b$ has zero average in $x \in \T$, as  deduced easily  from  (\ref{LSP}b) and the Dirichlet condition $b$. It follows that $b$ belongs to $L^2_{loc}(\R_+, H^1(\Omega))$.

These {\it a priori} estimates, combined with a classical approximation procedure, allow to state the following well-posedness result:
\begin{prop}
Assume that $U \in W^{2,\infty}(\R_+)$, $z U', z U''  \in L^\infty(\R_+)$. Let $u_0 \in L^2(\Omega)$ s.t.  $\omega_0 = \pa_z u_0 \in L^2(\Omega)$, $u_0\vert_{z=0} = 0$. Let $b_0 \in L^2_{loc}(\Omega)$ s.t. $\pa_z b_0 \in L^2(\Omega)$, $b_0\vert_{z=0} = 0$ and with zero average in $x$.  Then there there exists a unique solution $(u,v,b)$ of \eqref{LSP} satisfying \eqref{Gronwall1}-\eqref{Gronwall2}, 
$(u,b)\vert_{t=0} = (u_0, b_0)$. 
\end{prop}
\begin{rmk}
The main point of the proposition is that it does not involve any monotonicity assumption on the velocity profile $U$. This is in sharp contrast with the usual Prandtl system and its linearizations. In particular, when $U$ has a non-degenerate critical point, system \eqref{LP} does not admit this kind of solutions, see \cite{GeNg}. The difference comes from the control of $\pa_x u$ provided by the relation of Shercliff type. Let us stress that there is even a regularization effect in $x$, as no regularity in $x$ is required at initial time. 
\end{rmk}
\subsection{Fully nonlinear MHD layer}
In the specific regime in which $\Rey \sim \Rm \sim \Ha$, the formal model governing the boundary layer is \eqref{MHDBL}. Its 2D version reads 
\begin{equation} \label{2DMHDBL} 
\left\{
\begin{aligned}
& \pa_t u  + u \pa_x u + v \pa_z u  -   \pa^2_z u =   S \bb \cdot \na b  - \pa_x p^\infty, \\
&  \pa_t \bb   - \na^\perp (\bu \times \bb) - \frac{\Rey}{\Rm} \pa^2_z  \bb  = 0, \\
& \pa_x u + \pa_z v = \dive \bb = 0, \\
&   u\vert_{z=0} = v\vert_{z=0}, \quad  \bb\vert_{z=0} = \be_x, \\
 & u \rightarrow u^\infty, \quad b \rightarrow b^\infty, \quad \text{ as } z \rightarrow +\infty. 
\end{aligned}
\right. 
\end{equation}
We recall that $\bu = (u,v)$ and $\bb = (b,c)$ are the 2D velocity and magnetic fields respectively. We also recall that the cross product of $\bu$ and $\bb$ is a scalar function: $\bu \times \bb = uc - bv$. To investigate the stability properties of this system, we consider once more a simple linearization, around 
\begin{equation} \label{ref_solution}
 u = U(z), \quad v = 0, \quad \bb = \be_x.
 \end{equation}
The linearized equations are 
\begin{equation} \label{LMHDBL} 
\left\{
\begin{aligned}
& \pa_t u  + U \pa_x u + U' \, v   -   \pa^2_z u =   S  \pa_x  b, \\
&  \pa_t \bb   - \na^\perp (v - U \, c) - \frac{\Rey}{\Rm} \pa^2_z  \bb  = 0, \\
& \pa_x u + \pa_z v =  \dive \bb = 0, \\
&   u\vert_{z=0} = v\vert_{z=0}, \quad  \bb\vert_{z=0} = 0, \\
 & u \rightarrow 0, \quad b \rightarrow 0, \quad \text{ as } z \rightarrow +\infty. 
\end{aligned}
\right. 
\end{equation}
Here, $\bu = (u,v)$ and $\bb = (b,c)$ are the perturbations of the reference solution \eqref{ref_solution}. 

Note that by the conditions  $\pa_x b + \pa_z c = 0$, $c\vert_{z=0} = 0$, $c$ has zero average in $x$. Moreover, the evolution of the $x$-average of $b$ is decoupled and solves 
$$ \pa_t \int_\T b - \frac{\Rey}{\Rm} \pa^2_z \int_\T b = 0. $$
Hence, {\em there is no loss of generality in assuming that $b$ has zero average in $x$ as well}. With regards to the divergence-free condition, this means we can write $\bb = \na^\perp \phi$, for some function $\phi$ which is periodic with zero average in $x$.   We can then write the second component of  (\ref{LMHDBL}b)  as 
$$ \pa_t \pa_x \phi - \pa_x (v - U \, \pa_x \phi)  -  \frac{\Rey}{\Rm}  \pa^2_z \pa_x \phi = 0 $$
or equivalently 
\begin{equation} \label{streamphi}
\pa_t \phi + U \pa_x \phi - v -  \frac{\Rey}{\Rm}  \pa^2_z \phi = 0. 
\end{equation}
This last equation is a key ingredient in the stability analysis of \eqref{LMHDBL}. The idea is that, combining the equation (\ref{LMHDBL}a) on $u$ and \eqref{streamphi}, one can get rid of the bad term in $v$, responsible for the possible loss of one derivative in $x$. This idea is reminiscent of article \cite{MaWo} about the classical Prandtl equation. In \cite{MaWo}, a similar cancellation of the $v$ term was obtained combining the equations on $u$ and $\omega = \pa_y u $. In the linearized setting, the appropriate combination was 
$g = \omega - \frac{U''}{U'} u$. However, some monotonicity of the velocity profile  was needed, in order to divide by $U'$. {\em The main point in the present MHD context is that no monotonicity of the velocity profile is needed to obtain well-posedness}. We rather consider the following modified velocity: 
$$ \tilde u =   u + U' \phi. $$
Summing (\ref{LMHDBL}a)  and $U' \times$\eqref{streamphi}, we get 
\begin{equation} \label{eq_tildeu}
\pa_t \tilde u + U \pa_x \tilde u - \pa^2_z \tilde u =  S  \pa_x  b + \frac{\Rey}{\Rm}  U' \pa^2_z \phi - \pa^2_z \left( U' \phi \right),
 \end{equation}
while the equation on $b = \bb \cdot \be_x$ can be written as 
\begin{equation} \label{eq_b}
\pa_t b + U \pa_x b - \pa_x \tilde u  - \frac{\Rey}{\Rm} \pa^2_z b  = 0. 
\end{equation}
Formulation \eqref{eq_tildeu}-\eqref{eq_b} is much better behaved than the original formulation, and will allow to establish stability. Indeed,  a standard energy estimate yields 
$$\frac{d}{dt} \left(  \frac{1}{2}  \| \tilde u \|_{L^2}^2 + \frac{S}{2} \| b \|_{L^2}^2 \right)  \: + \: \| \pa_z \tilde u \|_{L^2}^2 +  \frac{S \Rey}{\Rm}  \| \pa_z b \|_{L^2}^2 \le  \frac{\Rey}{\Rm} \int_\Omega  U' (\pa^2_z \phi) \tilde u  - \int_\Omega  \pa^2_z \left( U' \phi \right) \, \tilde u,  $$
where we have used the identity 
$$-S\int_\Omega \pa_x \tilde u \, b =  S \int_{\Omega} \pa_x b \, \tilde u. $$
To control the r.h.s., we then use that $ \pa_z \phi = - b$. In particular, 
$$ \| \pa_z \phi \|_{L^2} = \| b \|_{L^2}, \quad  \| \pa^2_z \phi  \|_{L^2} = \| \pa_z b \|_{L^2}, \quad \| z^{-1} \phi \|_{L^2} \le 2 \| b \|_{L^2}.  $$
Hence, 
$$ \frac{d}{dt} \left(  \frac{1}{2}  \| \tilde u \|_{L^2}^2 + \frac{S}{2} \| b \|_{L^2}^2 \right)  \: + \: \| \pa_z \tilde u \|_{L^2}^2 +  \frac{S \Rey}{\Rm}  \| \pa_z b \|_{L^2}^2 \: \le \:  C (\| b \|_{L^2}  + \| \pa_z b \|_{L^2}) \| \tilde u \|_{L^2}, $$
where the constant $C$ depends implicitly on $\| U' \|_{L^\infty}$, $\| U'' \|_{L^\infty}$, $\| z U''' \|_{L^\infty}$. After application of Young's inequality: 
$$ \frac{d}{dt} \left(  \frac{1}{2}  \| \tilde u \|_{L^2}^2 + \frac{S}{2} \| b \|_{L^2}^2 \right)  \: + \: \| \pa_z \tilde u \|_{L^2}^2 +  \frac{S \Rey}{2\Rm}  \| \pa_z b \|_{L^2}^2 \: \le \:  C' \left(\| \tilde u \|_{L^2}^2 + \| b \|_{L^2}^2 \right) $$
Gronwall inequality yields 
$$    \| \tilde u(t) \|_{L^2}^2 +  \| b(t) \|_{L^2}^2 \: + \: \int_0^t \left(  \| \pa_z \tilde u \|_{L^2}^2 +   \| \pa_z b \|_{L^2}^2 \right)    \le M  \left(  \| \tilde u(0) \|_{L^2}^2 +  \| b(0) \|_{L^2}^2 \right) e^{Mt}, \quad \forall t \ge 0   $$  
where $M > 0$ is large enough.  Using $\| (U', U'') \phi \|_{L^2} \le 2 \| z (U',U'') \|_{L^\infty} \, \| b \|_{L^2}$, it follows that  
 \begin{equation} \label{Gronwall_MHDBL}
  \| u(t) \|_{L^2}^2 +  \| b(t) \|_{L^2}^2 \: + \: \int_0^t \left(  \| \pa_z  u \|_{L^2}^2 +   \| \pa_z b \|_{L^2}^2  \right)   \le 
   M' \left(  \| u(0) \|_{L^2}^2 +  \| b(0) \|_{L^2}^2 \right) e^{M' t} \quad \forall t \ge 0
 \end{equation}
for some  $M'$ large enough. 

As in the case of system \eqref{LSP}, we can combine the previous estimate with a standard approximation procedure, and obtain the well-posedness of \eqref{LMHDBL}:
\begin{prop}
Assume that $U \in W^{3,\infty}(\R_+)$, $z U', z U'', z U'''  \in L^\infty(\R_+)$.Let $u_0 \in L^2(\Omega)$. Let $\phi_0 \in L^2_{loc}(\Omega)$, such that $b_0  = \pa_y \phi_0 \in L^2(\Omega)$,
 $\phi_0\vert_{z=0} = 0$  and with zero average in $x$. Then there exists a unique solution of \eqref{LMHDBL} satisfying \eqref{Gronwall_MHDBL}, 
$u\vert_{t=0} = 0$, $\bb\vert_{t=0} = -\na^\perp \phi_0$.  
\end{prop}
\begin{rmk}
The velocity and magnetic vertical components $v$ and $c$  provided by this well-posedness  proposition have  weak regularity with respect to $x$. For instance, $v = -\int_0^y \pa_x u$ has to be understood as the $x$ derivative of a function in  $L^2(\T, H^2_{loc}(\R_+))$. For more regularity, one should impose more $x$ regularity on the data. 
\end{rmk}
\begin{rmk}
While completing the writing of this work, we got aware of the independent recent work \cite{Tong2} by Cheng-Jie Liu, Feng Xie and Tong Yang. These authors consider the same system as \eqref{2DMHDBL}, with the insulating boundary replaced by a conducting one, which amounts to replacing the condition $b\vert_{z=0} = 0$ by  $\pa_z b\vert_{z=0} = 0$. They establish well-posedness in Sobolev spaces for the nonlinear system, through a change of unknowns which is a nonlinear analogue of our $\tilde u$.   
\end{rmk}
\subsection{Conclusion}
We achieved a formal derivation and stability analysis of boundary layer models in MHD. This work was motivated by some contradictory results  on the stabilizing or destabilizing role of the magnetic field, notably when it is tangent to the boundary. The boundary layer models are in most regimes linear, but for some asymptotics  of the parameters, the role of the nonlinearities can not be ignored, leading to models of Prandtl type with extra magnetic features.  We investigated the stability to high frequencies of these nonlinear models, restricting to simple linearizations. Our analysis shows that in the case of tangent magnetic fields, the growth rate of high tangential frequencies is no longer growing with the wave number, contrary to what happens for the Prandtl system when the velocity has inflexion points. It favours the idea of  stabilization by  the magnetic field. 

\bibliographystyle{abbrv}

\end{document}